\documentclass[11pt, reqno]{amsart}
\usepackage{amscd}
\usepackage{amssymb}
\usepackage{amsmath}
\usepackage{enumerate}

\begin{document}
\title [Intrinsic operators ] {Intrinsic Operators from Holomorphic Function Spaces to Growth Spaces}
\author{Nina Zorboska }

\address{Department of Mathematics, University of Manitoba, Winnipeg, MB, R3T 2N2} 
\email{zorbosk@cc.umanitoba.ca}
\thanks{ Research supported in part by NSERC grant.}

\begin{abstract}
We determine the boundedness and compactness of a large class of operators, mapping from general Banach spaces of holomorphic functions into a particular type of spaces of functions determined by the growth of the functions, or the growth of the functions derivatives. The results show that the boundedness and compactness of such intrinsic operators depends only on the behaviour on the point evaluation functionals. They also generalize previous similar results about several specific classes of operators, such as the multiplication, composition and integral operators. 
\end{abstract}

\maketitle

\def\C{\mathbb{C}}
\def\D{\mathbb{D}}
\def\N{\mathbb{N}}
\def\R{\mathbb{R}}
\def\a{\alpha}
\def\b{\beta}
\def\g{\gamma}
\def\d{\delta}
\def\e{\varepsilon}
\def\l{\lambda}
\def\r{\rho}
\def\s{\sigma}
\def\t{\theta}
\def\B{\mathcal{B}}
\def\H{\mathcal{H}}
\def\p{\phi}
\def\vp{\varphi}

{\bf{1. Introduction}}

\bigskip

In this article we characterize the boundedness, compactness and weak compactness of a general class of operators acting on spaces of holomorphic functions, and mapping into the class of so called growth spaces. 
While there is an extensive existing literature on properties of operators acting on function spaces, it mostly deals with either specific types of operators, or with specific types of spaces. Few of the more recent relevant papers of this kind are, for example, [4], [5], [9], [10]. 
Our goal here is to consider the boundedness and compactness of a natural general class of operators on Banach spaces of holomorphic functions by using fairly general techniques, showing also that a large number of the existing specific results follow from our more general consideration.

Let $\D$ denote the open unit disk in the complex plane. Let $\H(\D)$ be the space of all holomorphic functions
on $\D$ with the topology of uniform convergence on compact subsets of $\D$, which we denote by $\tau_{uc}$.

We will say that an operator $T: \H(\D) \to \H(\D)$ is \textit{intrinsic} for $\H(\D)$ if it is well defined and continuous, namely when it maps 
$\tau_{uc}$ convergent sequences in $\H(\D)$ into $\tau_{uc}$ convergent sequences. There are many specific classes of intrinsic operators whose properties have been extensively studied in the last few decades. They include operators related naturally  to the operations that one can perform on a holomorphic function, such as multiplication, composition, differentiation or integration. Since as a consequence of Montel's theorem, the identity operator on $\H(\D)$ is compact, the really interesting questions arise when we consider intrinsic operators restricted to specific subspaces of $\H(\D)$. Our main goal is to determine the boundedness and compactness of these operators, when mapping a fairly general class of Banach spaces contained in $\H(\D)$ into the growth spaces $H_v$ 
and $\B_v$ defined below. 

\medskip

Recall that a \textit{weight} is a continuous function $v: \D \to (0, 1]$. A weight $v$ is is called a 
\textit{radial weight} if it is also a radial function, i.e. when $v(z)=v(|z|)$. If moreover $\lim_{|z|\to 1} v(z)=0$ and $v$ is non-increasing with respect to 
$|z|$, then the weight is usually referred to as a \textit{typical weight}. Unless stated otherwise, we will assume that our weights are typical weights.
The \textit{growth spaces} corresponding to a weight $v$ are defined as follows:
$$H_v= \{ f \in \H(\D); ||f||_{H_v}=\sup_{z \in \D} v(z) |f(z)| < \infty\},$$ 
$$\B_v = \{ f \in \H(\D); ||f||_{\B_v}=\sup_{z \in \D} v(z) |f'(z)| < \infty \}.$$

The growth spaces are Banach spaces of holomorphic functions; $||f||_{H_v}$ is a norm, while $||f||_{\B_v}$ is a semi-norm. They are both non-separable spaces, containing the polynomials, and the corresponding closure of polynomials are the so called \textit{little growth spaces}:
$$H_{v,0}= \{ f \in \H(\D); \lim_{|z| \to 1} v(z) |f(z)| =0 \},$$ 
$$\B_{v,0} = \{ f \in \H(\D); \lim_{|z| \to 1} v(z) |f'(z)|=0 \}.$$

If the growth function $v(z)= (1-|z|^2)^{\b}, \b>0$, the growth spaces $H_v$ and $\B_v$ are usually denoted by $H_{\b}$ and $\B_{\b}$ correspondingly, and we will refer to them as classical growth spaces. The spaces $\B_{\b}$ are also called Bloch-type spaces.
We also note that when  $\b>1$, the identity operator acts as an isomorphism between the space $H_{\b -1}$ and the space $\B_{\b}^0$, a subspace of $\B_{\b}$ of functions $f$ with $f(0)=0$. 

\bigskip

Growth spaces are an interesting and important class of Banach spaces of holomorphic functions. They have been explored in many different contexts and there are many general and more specific references such as, for example 
[2], [3], [14].
Some well known properties of these spaces, the proofs of which can be found in the cited papers and their further references, are that:

$\bullet$ For a radial weight $v$, $H_v$ is strictly bigger than $H^{\infty}$ (the space of bounded holomorphic functions on $\D$) if and only if 
$\lim_{|z|\to 1} v(z)=0$. If $\limsup_{|z|\to 1} v(z)>0$, then $H_{v, 0}=\{0\}$.

$\bullet$ The topologies on $H_v$ and $\B_v$ are stronger than the $\tau_{uc}$ topology;

$\bullet$  The double duals $H_{v,0}^{**}$ and $\B_{v,0}^{**}$ are isometrically isomorphic to $H_v$ and $\B_v$ correspondingly;

$\bullet$ The point evaluation functionals on $H_{v,0}$ and $\B_{v,0}$ are bounded and are uniquely extended to point evaluation functionals on $H_{v}$ and $\B_{v}$ with equal norms. We will denote them correspondingly by $K_z^H$ and $K_z^{\B}$.

$\bullet$ The differentiation operator sending $f \to f'$ is an isometric isomorphism between $\B_v^0$, the subspace of $\B_v$ of functions with $f(0)=0$, and $H_v$.

$\bullet$ The differentiation operator sending $f \to f'$ is an isomorphism between $H_{\b}^0$, the subspace of 
$H_{\b}$ of functions with $f(0)=0$, and $H_{\b+1}$.

$\bullet$ The derivative point evaluations $K_{z,1}^{\B}$ on $\B_{v, 0}$ (and on $\B_v$) are bounded linear functionals.

$\bullet$ The maps $z \to K_z^H$, $z \to K_z^{\B}$ and $z \to K_{z,1}^{\B}$ are continuous, and the norms of the point evaluations and derivative point evaluations go to infinity, as $|z| \to 1$.

\bigskip

The domain for the operators considered will be of a fairly general type, and will be referred to as an initial space. More precisely, we will say that $X$ is an \textit{initial space} if 
$X \subset \H(\D)$ is a Banach space containing the polynomials, and such that the closed unit ball 
$B_X$ of $X$ is compact with respect to the $\tau_{uc}$ topology, namely: 

$\bullet$ for every sequence $\{f_n\}$ in the closed unit ball $B_X$ of $X$, there exists a subsequence $\{f_{n_k}\}$ and a function $f$ in $X$ such that $f_n \to f$ in the $\tau_{uc}$ topology.

This includes a large class of well known classical spaces of analytic functions such as the Hardy, weighted Bergman, weighted Dirichlet, Besov-type, BMOA and the growth spaces. The little growth spaces on the other hand are not initial spaces since the closed unit ball of the little growth spaces is only relatively compact with respect to the 
$\tau_{uc}$ topology, namely it is $\tau_{uc}$ dense in the closed unit ball of the corresponding growth spaces (see [3] and [9]).

\smallskip

As we will show below, the boundedness and the compactness of intrinsic operators is completely determined by their behaviour on the point evaluation functions (also referred to as kernel functions). Section 2 contains the general boundedness criteria and several specific examples of applications of these criteria. Some of the applications to specific classes of operators give results that have been known before, and some of the applications are new results. This section also contains a result about a special class of non-reflexive spaces, for which the boundedness and weak compactness of the intrinsic operator on the space and on its double dual are closely related. The general compactness and weak compactness characterization of the intrinsic operators, and few further specific examples and applications of these results are presented in Section 3.

\medskip

{\bf{2. Boundedness}}

\bigskip

Recall that $K_z^H$ denotes a bounded point evaluation functional for $H_v$, namely 
$K_z^H(g) = g(z), \forall g \in H_v.$
Thus, for any $z$ in $\D$ 
$$||K_z^H|| = \sup_{z \in \D}\{|K_z^H(g)|; ||g||_{H_v} \le 1\}= \sup_{z \in \D}\{|g(z)|; ||g||_{H_v} \le 1\} \le \frac{1}{v(z)}.$$ 

Similarly, it is also easy to see that the norm of the derivative point evaluation functional $K_{z, 1}^{\B}$ for the space $\B_v$ is such that 
$||K_{z, 1}^{\B}|| \le \frac{1}{v(z)}.$

\medskip

In order to give a boundedness criteria for general linear transformations on subspaces of $\H(\D)$, we introduce the following notation. For a linear transformation
$T: \H(\D) \to \H(\D)$, a Banach space $X \subset \H(\D)$ and $z \in \D$, define a linear map from $X$ into $\C$ by:
$$(K_z^H)_T (f) = Tf(z).$$
Hence, $(K_z^H)_T \in X^{\dagger}$, where $X^{\dagger}$ denotes the algebraic dual of $X$. Similarly, let $(K_{z, 1}^{\B})_T$ be the linear map in $X^{\dagger}$ defined by 
$$(K_{z, 1}^{\B})_T(f)= (Tf)^{\prime}(z).$$

We have the following operator boundedness characterizations when the target space is a growth space. The statements and the proofs of the next two main theorems in this section are mostly just a general characterization of boundedness of linear operators on Banach spaces, using algebraic duals and operator adjoints. We provide the statements and the proofs in this specific context in order to show more clearly the notation choices, and also to show how the existing boundedness characterizations for some specific operators follow as a consequence of this more general case.

\bigskip

\textbf{Theorem 2.1.} \textit{Let $T: \H(\D) \to \H(\D)$ be a linear transformation, $X \subset \H(\D)$ a Banach space, $v$ a radial weight function, and let
$(K_z^H)_T$ and $(K_{z, 1}^{\B})_T$ be as defined above. Then}

 \textit{(i) $T: X \to H_v$ is bounded if and only if 
\begin{equation*}
(K_z^H)_T \in X^*, \forall z \in \D  \quad \mbox{and} \quad
\sup_{z \in \D}v(z) ||(K_z^H)_T|| < \infty.
\tag{2.1.1}
\end{equation*}}

\textit{(ii) $T: X \to \B_v$ is bounded if and only if 
\begin{equation*}
(K_0^{\B})_T \in X^*, (K_{z,1}^{\B})_T \in X^*, \forall z \in \D  \quad \mbox{and} \quad
\sup_{z \in \D}v(z) ||(K_{z,1}^{\B})_T|| < \infty.
\tag{2.1.2}
\end{equation*}}

\begin{proof} 

\textit{(i)} Let $T: X \to H_v$ be bounded. Then $T^*: H_v^* \to X^*$ is also bounded, and if  $C=||T^*||=||T||$, then
$||T^* K_z^H|| \leq C ||K_z^H||, \forall z \in \D$. Also, since $Tf \in H_v, \forall f \in X$, we have that  for all $f$ in $X$
$$(K_z^H)_T (f)=Tf(z)=K_z^H(Tf)= T^*K_z^H(f),$$ 
and so $(K_z^H)_T= T^*K_z^H \in X^*, \forall z \in \D$.
Moreover, since $||K_z^H|| \leq 1/v(z), \forall z \in \D$, 
$$\sup_{z \in \D}v(z)||(K_z^H)_T||= \sup_{z \in \D}v(z) ||T^{*} K_z^H|| \leq C \sup_{z \in \D}v(z)||K_z^H|| \leq C< \infty.$$
For the other direction, assume that (2.1.1) holds, i.e. that $(K_z^H)_T \in X^*, \forall z \in \D$ and that $\exists C>0$ such that
$||(K_z^H)_T|| \leq C/v(z), \forall z \in \D$. 
But then for all $f \in X, z \in \D$ we have that 
$$|Tf(z)|= |(K_z^H)_T(f)| \leq ||(K_z^H)_T|| ||f||_X \leq C \frac{1}{v(z)} ||f||_X,$$ 
and so (2.1.1) implies that 
$$||Tf||_{H_v}= \sup_{z \in \D}v(z) |Tf(z)| = \sup_{z \in \D}v(z) |(K_z^H)_T(f)| \leq C ||f||_X,$$
i.e. $T: X \to H_v$ is bounded.

\textit{(ii)} The proof of the second part is similar, and we leave it to the reader. We note only that in this case 
$T: X \to \B_v$ bounded implies
that $(K_{z,1}^{\B})_T=T^*K_{z,1}^{\B}$, and that in general, for $f \in X$, $||Tf||_{\B_v}=\sup_{z \in \D}v(z) |(K_{z,1}^{\B})_T(f)|$.
\end{proof}
\medskip

Note that from the inequalities in the proof of part $(i)$ above, and since for $f \in X, ||f||_X \leq 1$, we have that 
$v(z)|T^*K_Z^H(f)| \leq v(z)||T^*K_z^H||$ for all $z \in \D$, we can easily see by
taking appropriate supremums that when $T:X \to H_v$ is bounded, then 
$$||T||=\sup_{z \in \D}v(z) ||(K_z^H)_T||.$$

We get a similar conclusion also for part $(ii)$ of the previous theorem. In this case though, we are dealing with semi-norms, and so we need to restrict the range of $T$ to a  subspace of $\B_v$. Namely, if $P$ is the map onto the subspace $\B_v^0$ defined by $Pf=f-f(0)$ and if $T:X \to \B_v$ is bounded, then $||PT||=\sup_{z \in \D}v(z) ||(K_{z,1}^{\B})_T||$.

\bigskip

We will illustrate the application of the previous theorem first to two specific classes of operators, i.e. weighted composition operators and integral operators. There is a vast literature exploring the properties of these two classes of operators, acting on various spaces of holomorphic functions. 
One of the most recent articles of this kind is, for example, [4], and we cite below some of its results which also follow from Theorem 2.1.
Another recent relevant interesting article is [5]. It characterizes the bounded, weakly compact and compact integral operators from a variety of specific Banach spaces of holomorphic functions into the space of bounded holomorphic functions $H^{\infty}(\D)$, and into the disk algebra $A(\D)$. Even though its general approach and the actual target spaces are different, few of the particular ideas were an inspiration for the generalizations obtained here.

\medskip

Let $u \in \H(\D)$ and let $\phi:\D \to \D$ be holomorphic. Define a weighted composition operator $W_{u,\phi}$ on $X \subset \H(\D)$ by
$$W_{u,\phi}f=uf\circ\phi, \forall f \in X.$$
The operator $W_{u,\phi}$ is an intrinsic operator on $\H(\D)$. 

If $X$ is also a Banach space with bounded point evaluations, i.e.
$$\forall z \in \D, \exists \delta_z \in X^{*} ~\mbox {such that}~ \delta_z(f)=f(z), \forall f \in X,$$
then the condition $(K_z^H)_T \in X^*$, for $T$ a weighted composition operator, is equivalent to $u(z) \delta_{\phi(z)} \in X^*$,  since for any $f \in X$
$$(K_z^H)_T(f)=Tf(z)=W_{u,\phi}f(z)=u(z)f(\phi(z))=u(z)\delta_{\phi(z)}(f).$$
For any $z \in \D$, $u(z) \delta_{\phi(z)}$ is trivially in $X^*$, since $\phi(z) \in \D$ and $X$ is a Banach space with bounded point evaluations.
Furthermore, the previous equalities also show that  $W_{u,\phi}^* K_z^H=u(z)\delta_{\phi(z)}$. 

The next example is the example of an integral operator. 
Let $g \in \H(\D)$ and let $X \subset \H(\D)$. The integral operator $T_g$ is defined by
$$T_g f(z)= \int_0^z f(w)g'(w)dw, \forall f \in X, z \in \D.$$
It is easy to see that the integral operator $T_g$ is also an intrinsic operator on $\H(\D)$. Note also that 
$T_g f(0)= 0, \forall f \in X$.

The differential operator $Df=f'$ is an isometric isomorphism between the space 
$\B_v^0=\{f \in B_v; f(0)=0\}$ and the space $H_v$. Thus, it is easy to see that the integral operator $T_g: X \to \B_v$ is bounded if and only if the multiplication operator $M_{g'}: X \to \H_v$ is bounded. 
Since for $\b>0$ the differential operator is an isomorphism between the space $H_{\b}^0$  and $H_{\b+1}$, we get similarly that 
$T_g: X \to H_{\b}$ is bounded if and only if  $M_{g'}:X \to H_{\b+1}$ is bounded.
Also, from the discussion above on weighted composition operators, we have that if $X \subset \H(\D)$ is a Banach space with bounded point evaluations $\delta_z$, then the condition $(K_z^H)_{M_{g'}} \in X^*$ is equivalent to $g' (z) \delta_{z} \in X^*$, which is trivially true, and that $M_{g'}^* K_z^H = g'(z) \delta_z$. 

\medskip

Using these few general facts about the weighted composition operators and the integral operators, the next three corollaries follow as an application of the operator boundedness characterization from Theorem 2.1.

\bigskip

\textbf{Corollary 2.1.[4]} \textit{Let $v$ be a radial weight function, let $X \subset \H(\D)$ be a Banach space with bounded point evaluations $\delta_z$, and let $\b>0$. Then:} 

\textit{(i) The weighted composition operator $W_{u,\phi}: X \to H_v$ is bounded 
if and only if 
$\sup_{z \in \D}v(z) |u(z)| ||\delta_{\phi(z)}|| < \infty$. When  $W_{u,\phi}: X \to H_v$ is bounded, then
$$ ||W_{u,\phi}|| = \sup_{z \in \D}v(z) |u(z)| ||\delta_{\phi(z)}||.$$ }

\textit{(ii) The integral operator $T_g: X \to \B_v$ is bounded if and only if 
\linebreak
$\sup_{z \in \D}v(z) |g'(z)| ||\delta_z|| < \infty$. When
$T_{g}: X \to \B_v$ is bounded, then
$$ ||T_g|| =  \sup_{z \in \D}v(z) |g'(z)| ||\delta_z||.$$ }

\textit{(iii) The integral operator $T_g: X \to H_{\b}$ is bounded if and only if 
\linebreak
$\sup_{z \in \D}(1-|z|^2)^{\b+1} |g'(z)| ||\delta_z|| < \infty$. When
$T_{g}: X \to H_{\b}$ is bounded, then
$$ ||T_g|| =  \sup_{z \in \D}(1-|z|^2)^{\b+1} |g'(z)| ||\delta_z||.$$}

Applying Corollary 2.1 to particular choices of Banach spaces with bounded point evaluations gives further specific characterizations of bounded weighted composition and integral operators. One only needs to know the norm of the point evaluation functionals for the corresponding space $X$.
We illustrate this in the next corollary, by stating some of the concrete results for the integral operator acting on the Hardy and weighted Bergman spaces. We also include few further specific comments, added in order to show the commonalities with some of the results for the integral operators mapping into the space $H^{\infty}$ (see parts of Theorem 1.3 and Theorem 2.6 in [5]).

For more details and more examples of this kind for the weighted composition operators or for the integral operators with $X$ being either the space $BMOA$, $VMOA$, the Besov,  or  the Bloch spaces, see further [4] and [5].

\medskip   

Recall first that for $1 \le p < \infty$, the space $H^p$ is the classical Hardy space of functions in $\H(\D)$, with bounded point evaluations $\delta_z$  such  that $||\delta_z||=(1-|z|^2)^{-\frac{1}{p}}$.

For $\a>-1$ and $0<p<\infty$, $A_{\a}^p$ denotes the weighted Bergman space of functions in $\H(\D)$. Namely, $A_{\a}^p$ is a 
$L^p(d_{\mu_{\a}})$ Banach space with $d_{\mu_{\a}}(z)=(1-|z|^2)^{\a}dA(z)$, where $dA$ refers to the normalized Lebesque area measure on $\D$. Furthermore, recall that $A_{\a}^p$ is a Banach space of bounded point evaluations with 
$$||\delta(z)|| = (1-|z|^2)^{-\frac{\a+2}{p}}.$$

An application of the previous corollary and the comments above give the following more specific boundedness characterizations for the integral operators. For similar characterizations for weighted composition operators on the Hardy and weighted Bergman spaces see [4].

\bigskip

\textbf{Corollary 2.2.} \textit{Let $g \in \H(\D)$, $1< p <\infty$, $\a>-1$ and $\b>0$.}

\textit{(i) The integral operator $T_g: H^p \to H_{\b}$ is bounded if and only if  $g \in \B_{\beta+1 - \frac{1}{p}}$, i.e
$$\sup_{z \in \D}(1-|z|^2)^{\beta+1 - \frac{1}{p}} |g'(z)| < \infty.$$}

\textit{(ii) The integral operator $T_g: H^p \to \B_{\b}$ is bounded if and only if 
$$\sup_{z \in \D}(1-|z|^2)^{\beta - \frac{1}{p}} |g'(z)| < \infty.$$
Thus, if $p > \frac{1}{\beta}$ the operator $T_g: H^p \to \B_{\b}$ is bounded if and only if $g \in \B_{\beta - \frac{1}{p}}$,\linebreak
When $p =\frac{1}{\beta}$, then $T_g: H^p \to \B_{\b}$ is bounded if and only if $g$ is a Lipschitz function, and
if $p < \frac{1}{\beta}$, then $T_g: H^p \to \B_{\b}$ is bounded only when $g$ is a constant.}

\textit{(iii) The integral operator $T_g: A_{\a}^p \to H_{\b}$ is bounded if and only if 
$$\sup_{z \in \D} (1-|z|^2)^{\beta+1 - \frac{2+\a}{p}} |g'(z)| < \infty.$$
Thus, if $p >\frac{\a+2}{\b+1}$ the operator $T_g: A_{\a}^p \to H_{\b}$ is bounded if and only if $g \in \B_{\beta+1 - \frac{2+\a}{p}}$. 
When $p =\frac{\a+2}{\b+1}$, then $T_g: A_{\a}^p \to H_{\b}$ is bounded if and only if $g$ is a Lipschitz function, and
if $p <\frac{\a+2}{\b+1}$, then $T_g: A_{\a}^p \to H_{\b}$ is bounded only when $g$ is a constant.}

\textit{(iv) The integral operator $T_g: A_{\a}^p \to \B_{\b}$ is bounded if and only if 
$$\sup_{z \in \D} (1-|z|^2)^{\beta - \frac{2+\a}{p}} |g'(z)| < \infty.$$
Thus, if $p >\frac{\a+2}{\b}$ the operator $T_g: A_{\a}^p \to \B_{\b}$ is bounded if and only if $g \in \B_{\beta - \frac{2+\a}{p}}$.
 When $p =\frac{\a+2}{\b}$, then $T_g: A_{\a}^p \to \B_{\b}$ is bounded if and only if $g$ is a Lipschitz function, and
if $p <\frac{\a+2}{\b}$, then $T_g: A_{\a}^p \to \B_{\b}$ is bounded only when $g$ is a constant.}

\bigskip

We mention one more interesting example of an integral-type operator, i.e. the generalized Ces\`aro operator, and state the specific boundedness results that follow when choosing the initial space $X$ to be a classical growth space. These results are related to the results in [1, Theorem 4.1] for the generalized Ces\`aro operator, and to results in [11] for the integral operator.

\medskip

For $g \in \H(\D)$, the generalized Ces\`aro operator $C_g$ is defined on the Banach space $X \subset \H(\D)$ by
$$C_g f(z)= \frac{1}{z} \int_0^z f(w)g'(w)dw, z \ne 0, ~\mbox{and}~ C_g f(0)=0.$$
When $g(z)= \log\frac{1}{1-z}$, i.e. $g'(z) = \frac{1}{1-z}$, $C_g$ is the classical Ces\`aro operator. 

The generalized Ces\`aro operator $C_g$ is closely related to the integral operator $T_g$. Namely, if $Y  \subset \H(\D)$ is such that the unilateral shift $S$ (i.e. the operator of multiplication by the identity function) is bounded on $Y$ and is left invertible, then
$T_g=S C_g$ and $C_g= \hat{S} T_g$, where $\hat{S}$ is the left inverse of $S$. Hence, $C_g:X \to Y$ is bounded on if and only if 
$T_g:X \to Y$ is bounded. Such are, for example, the classical growth spaces $H_{\b}$ and $\B_{\b}$.
 
One of the results in [1] determines the boundedness and compactness of generalized Ces\`aro operators with essentially rational symbols, acting on spaces satisfying the condition above, and which are further also contained in $H_{\b}$ for some $\b>0$ (see [1,Theorem 4.1]). The result in the next corollary generalizes this for the case when $X$ and $Y$ are classical growth spaces of the same kind.

\medskip

The corollary below follows from Theorem 2.1, the discussion before Corollary 2.1, and Corollary 2.1 itself. Recall also that for the growth space $H_{\b}$, the point evaluations norm is given by $||\delta_z||= (1-|z|^2)^{-\b}$, while for the Bloch-type spaces $\B_{\b}$, we have that 
\[ ||\delta_z|| \approx \begin{cases} 1, & 0<\b<1\\
\log\frac{1}{1-|z|^2}, & \b=1\\
\frac{1}{(1-|z|^2)^{\b-1}}, & \b>1.
\end{cases}\]
In the following result  $\mbox{Log}\B$ denotes the logarithmic Bloch space, i.e. the space of functions $f \in \H(\D)$ such that
$\sup_{z \in \D} (1-|z|^2) \log \frac{1}{1-|z|^2} |f'(z)| < \infty$.

\bigskip

\textbf{Corollary 2.3.} \textit{Let $g \in \H(\D)$ and let $\b, \gamma>0$.}

\begin{enumerate}[\textit{(a)}]
\item  \textit{The following are equivalent:}
 \begin{enumerate}
  
\item[\textit{(i)}]  \textit{$C_g: H_{\gamma} \to H_{\b}$ is bounded.}

\item[\textit{(ii)}]   \textit{$T_g: H_{\gamma} \to H_{\b}$ is bounded.}

\item[\textit{(iii)}]   \textit{$M_{g'}: H_{\gamma} \to H_{\b+1}$ is bounded.}

\item[\textit{(iv)}]   \textit{$\sup_{z \in \D} (1-|z|^2)^{\beta +1-\gamma} |g'(z)| < \infty.$}

\end{enumerate}
\textit{In the case when $\gamma=\b$, the condition (iv) is equivalent to $g \in \B$, where $\B=\B_1$ is the classical Bloch space.}

\item [\textit{(b)}] \textit{Let $\delta_z$ denote the point evaluation functional for $\B_{\gamma}$. Then the following are equivalent:}
 \begin{enumerate}
 
\item[\textit{(i)}]   \textit{$C_g: \B_{\gamma} \to \B_{\b}$ is bounded.}

\item[\textit{(ii)}]   \textit{$T_g: \B_{\gamma} \to \B_{\b}$ is bounded.}

\item[\textit{(iii)}]   \textit{$M_{g'}: \B_{\gamma} \to H_{\b}$ is bounded.}

\item[\textit{(iv)}]   \textit{$\sup_{z \in \D} (1-|z|^2)^{\beta}||\delta_z|| |g'(z)| < \infty.$}
\end{enumerate}
\textit{In particular, in the case when $\gamma=\b$, the condition (iv) is equivalent to:} 
\[ g \in \begin{cases} \B_{\b}, & 0<\b<1\\ 
\mbox{Log}\B, & \b=1\\
\B, & \b>1
\end{cases}\]
\end{enumerate}

\bigskip

The last comments in parts \textit{(a)} and \textit{(b)} above, for the case $\gamma=\b$ and for the generalized Ces\`aro operator (in the case of the unit ball in $\C^n$), appear as results in [6] and [11] correspondingly.  
Note also that the essentially rational functions considered as symbols for the generalized Ces\`aro operator in [1] are functions that belong to the Bloch space.

\bigskip

The next result gives a general boundedness characterization of operators mapping into the little growth spaces with typical weights.

\bigskip

\textbf{Theorem 2.2.} \textit{ Let $T: \H(\D) \to \H(\D)$ be a linear transformation , $X \subset \H(\D)$ a Banach space, and let $(K_z^H)_T$ and $(K_{z, 1}^{\B})_T$ be as defined above. Then}

 \textit{(i) $T: X \to H_{v, 0}$ is bounded if and only if 
\begin{equation*}
(K_z^H)_T \in X^*, \forall z \in \D  \quad \mbox{and} \quad
v(z) (K_z^H)_T \to 0 ~ \mbox{weak*, as} ~ |z| \to 1.
\tag{2.2.1}
\end{equation*}}

\textit{(ii) $T: X \to \B_{v,0}$ is bounded if and only if 
\begin{equation*}
(K_0^{\B})_T \in X^*, (K_{z,1}^{\B})_T \in X^*, \forall z \in \D  \quad \mbox{and} \quad
v(z) (K_{z,1}^{\B})_T \to 0 ~ \mbox{weak*, as} ~ |z| \to 1.
\tag{2.2.2}
\end{equation*}}

\begin{proof}

\textit{(i)} If $T: X \to H_{v, 0}$ is bounded, then $T^*: H_{v,0}^* \to X^*$ is also bounded, $(K_z^H)_T= T^*K_z^H$ and 
so $(K_z^H)_T \in X^*, \forall z \in \D$. Since $\forall f \in X, Tf \in H_{v,0}$, we have that 
$$v(z) |(K_z^H)_T(f)|=v(z)|Tf(z)| \to 0, ~\mbox{as} ~|z| \to 1,$$ 
i.e. $v(z) (K_z^H)_T$ converges weak$^*$ to $0$ in $X^*$.

For the other direction, using the closed graph theorem, it is enough to show that (2.2.2) implies that $Tf \in H_{v,0}, \forall f \in X$.
This follows directly, since 
$$\lim_{|z| \to 1} v(z) |Tf(z)|= \lim_{|z| \to 1}v(z)|(K_{z}^H)_T(f)|=0.$$
 
\textit{(ii)} The proof of the second part follows similarly and we leave it to the reader, recalling that 
$(K_{z,1}^{\B})_T(f)= (Tf)'(z), \forall f \in X$.

\end {proof}

\medskip

Recall that every weak$^*$ convergent sequence is bounded. Hence, every operator $T$ that maps $X$ into the little growth spaces is also bounded as an operator from $X$ into the corresponding growth spaces $H_v$ (or $\B_v$), as it should trivially be the case by the Closed Graph Theorem.

\bigskip

Following from Theorem 2.2 are the corresponding versions of corollaries 2.1, 2.2, and 2.3 for specific choices of operators and spaces, this time with the target space being replaced by the little growth spaces. We state here only the corresponding version of Corollary 2.1 for the weighted composition and integral operators on  general Banach spaces $X$ with bounded point evaluations, and leave the statements of the other specific corollaries for the reader. 

\bigskip

\textbf{Corollary 2.4.} \textit{Let $X \subset \H(\D)$ be a Banach space with bounded point evaluations $\delta_z$, and let $\b>0$. Then:} 

\textit{(i) The weighted composition operator $W_{u,\phi}: X \to H_{v,0}$ is bounded 
if and only if }
\linebreak
\hspace*{0.9cm} 
\textit{$v(z) u(z) \delta_{\phi(z)} \to 0 ~\mbox{weak}^*, \mbox{as} ~|z| \to 1.$} 

\textit{(ii) The integral operator $T_g: X \to \B_{v,0}$ is bounded if and only if} 
\linebreak
\hspace*{0.9cm} 
\textit{$v(z) g'(z) \delta_z \to 0 ~\mbox{weak}^*, \mbox{as} ~|z| \to 1.$}

\textit{(iii) The integral operator $T_g: X \to H_{\b,0}$ is bounded if and only if} 
\linebreak
\hspace*{0.9cm} 
\textit{$(1-|z|^2)^{\b+1} g'(z)\delta_z \to 0 ~\mbox{weak}^*, \mbox{as} ~|z| \to 1.$}

\bigskip

Note that since $X$ contains the constant function $1$, the weak$^*$ condition in the corollary above implies that: 
in case (i), the function $u$ belongs to $H_{v,0}$, and in the cases (ii), or (iii), that the function $g$ belongs to $\B_{v,0}$, or $H_{\b, 0}$ correspondingly. These are of course necessary conditions for the boundedness of the corresponding operators.

\bigskip

The general weak compactness criteria for the intrinsic operators will be addressed in the next section.
However, let us mention here that in case $X$ is a reflexive Banach space, every bounded operator $T$ on $ X$ is also weakly compact, and the above conditions also characterize the weakly compact intrinsic operators.
  
It is interesting though that for a special type of non-reflexive spaces $X$, with properties as specified below, the boundedness of the intrinsic operator on the double dual is closely tied to its weak compactness.
As it was shown in [5, Lemma 3.1] such is, for example, the  following special class of Banach spaces $X$: 

$X^*$ is separable, and $X^{**} \subset \H(\D)$ is such that:

(i) Every bounded sequence of functions in $X^{**}$ is uniformly bounded on compact subsets of $\D$, and

(ii) For every $z \in \D$, there exists $\delta_z \in X^*$ such that $F(z)=<F, \delta_z>$, for all $F$ in $X^{**}$.

Thus, the class of spaces $X$ in the next theorem includes, for example, the non-reflexive spaces $VMOA, H_{v,0}, \B_{v,0}$ and 
$A(\D)$, and the pre-duals of $H^{\infty}, BMOA, H_v$ and $\B_v$ (see [2], [5]).

\medskip

\textbf{Theorem 2.3.} \textit{Let $X$ be a Banach space such that the weak$^*$ topology on the unit ball of $X^{**}$ is the $\tau_{uc}$ topology, and let $T$ be an intrinsic operator on $\H(\D)$ such that 
$T: X \to H_{v,0} (\mbox{or}~ \B_{v,0})$ is bounded. Then, for all 
$F \in X^{**}$, $T^{**} F=T F$, and so $T: X^{**} \to H_v (\mbox{or}~ \B_{v})$ is bounded. Furthermore, the following are equivalent:}

\textit{(i) $T: X \to H_{v,0} (\mbox{or}~ \B_{v,0})$ is weakly compact.}

\textit{(ii) $T: X \to H_v (\mbox{or}~ \B_{v})$ is weakly compact.} 

\textit{(iii) $T: X^{**} \to H_{v,0} (\mbox{or}~ \B_{v,0})$ is bounded.}

\textit{(iv) $T: X^{**} \to H_{v} (\mbox{or}~ \B_{v})$ is weakly compact.}

\begin{proof} We will show the proof only for the case of $H_v$ and $H_{v,0}$. The proof for the case $\B_v$ and $\B_{v,0}$ is similar.

Since $T: X \to H_{v,0}$ is bounded and $H_{v,0}^{**}=H_v$, we have that $T^{**}: X^{**} \to H_{v}$ is weak$^*$ to 
weak$^*$ continuous (see [13, p. 29]). 

By Goldstine's theorem [13, p. 31], the closed unit ball $B_X$ of $X$ is weak$^*$ dense in
the closed unit ball $B_{X^{**}}$ of $X^{**}$. We have assumed that the weak$^*$ topology on the unit ball of $X^{**}$ is the $\tau_{uc}$ topology. Thus,
for any $F \in X^{**}$, there exists a sequence $\{f_n\}$ in $B_X$ such that $f_n \to F$ in $\tau_{uc}$, and 
since $T$ is intrinsic on $\H(\D)$, we have that $Tf_n= T^{**}f_n \to TF$ in $\tau_{uc}$, and so, also pointwise.
 
But since $T^{**}: X^{**} \to H_{v}$ is also weak$^*$ to weak$^*$ continuous, $T^{**}f_n \to T^{**}F$ weak$^*$ in $H_v$. The pointwise convergence in $H_v$ is weaker than the weak$^*$ convergence (see [3, p.103]), and so $T^{**}f_n \to T^{**}F$ pointwise. Thus, for every $z \in \D$, $T^{**}F(z)= TF(z)$, i.e. $T^{**} F=T F$.
 
Since $T: X \to H_{v,0}$ is bounded $T^{**}/X^{**}= T/X^{**}$ and $H_{v,0}^{**}=H_v$, we get that $T: X^{**} \to H_v$ is also bounded.

We have that (i) and (ii) are equivalent since $T(X) \subset H_{v,0}$ and $H_{v,0}$ is a closed subset of $H_v$. 

The equivalence of (i) and (iii)  follows from Gantmacher and Nakamura's theorem (see [8, p. 341]), since $T: X \to H_{v,0}$ is bounded,
$T^{**}/X^{**}= T/X^{**}$, and so $T: X \to H_{v,0}$ is weakly compact if and only if $T(X^{**}) \subset H_{v,0}$.

That (i) and (iv) are equivalent follows by Gantmacher's theorem, i.e. since $T: X \to H_{v,0}$ is bounded and $T^{**}/X^{**}= T/X^{**}$, 
$T: X \to H_{v,0}$ is weakly compact if and only if $T: X^{**} \to H_{v}$ is weakly compact.

\end{proof}

\medskip

As an illustration of the application of the previous theorem, we state a corollary with some results for three special choices for the space $X$. Recall that $VMOA^{**}=BMOA$, and that for a typical weight $w$, $H_{w,0}^{**}=H_w$ and $\B_{w,0}^{**}=\B_w$.

\bigskip

\textbf{Corollary 2.5} \textit{ Let $T$ be an intrinsic operator on $\H(\D)$ and let v and w be typical weights.}

\textit{(a) If $T: VMOA \to H_{v,0} (\mbox{or}~ \B_{v,0})$ is bounded, then}

\quad \textit{(i) $T: BMOA \to H_v (\mbox{or}~ \B_{v})$ is bounded, and}

\quad \textit{(ii) $T: BMOA \to H_{v,0} (\mbox{or}~ \B_{v,0})$ is bounded if and only if} 

\hspace{0.9cm} $T: BMOA \to H_{v,0} (\mbox{or}~ \B_{v,0})$ \textit{is weakly compact.} 

\textit{(b) If $T: H_{w,0} \to H_{v,0} (\mbox{or}~ \B_{v,0})$ is bounded, then}

\quad \textit{(i) $T: H_w \to H_v (\mbox{or}~ \B_{v})$ is bounded, and}

\quad \textit{(ii) $T: H_w \to H_{v,0} (\mbox{or}~ \B_{v,0})$ is bounded if and only if} 

\hspace{0.9cm} $T: H_w \to H_{v,0} (\mbox{or}~ \B_{v,0})$ \textit{is weakly compact.} 

\textit{(c) If $T: \B_{w,0} \to \H_{v,0} (\mbox{or}~ \B_{v,0})$ is bounded, then}

\quad \textit{(i) $T: \B_w \to H_v (\mbox{or}~ \B_{v})$ is bounded, and}

\quad \textit{(ii) $T: \B_w \to H_{v,0} (\mbox{or}~ \B_{v,0})$ is bounded if and only if} 

\hspace{0.9cm} $T: \B_w \to H_{v,0} (\mbox{or}~ \B_{v,0})$ \textit{is weakly compact.} 

\bigskip

{\bf{3. Compactness}}

\bigskip

In this section we explore the compactness and weak compactness of intrinsic operators mapping into the growth spaces. 

We start with the following general lemma which is a standard known result that shows up in the literature in several similar versions. The proof of the compactness can be found, for example, in [12, Lemma 3.7]. The proof of the weak compactness follows similarly, after using the Eberlein-$\check{\text S}$mulian theorem (see [5]). 

\bigskip

\textbf{Lemma 3.1.} \textit{Let $X$ be an initial space and let $Y \subset \H(\D)$ be a Banach space with bounded point evaluations. Let $T$ be an intrinsic operator on $\H(\D)$ such that $T: X \to Y$ is bounded. Then $T: X \to Y$ is compact (or weakly compact) if and only if for every bounded sequence $\{f_n\}$ in $X$ that converges to zero in $\tau_{uc}$, we have that $\{Tf_n\}$ converges in norm (or weakly) to zero.}

\bigskip

We have the following general characterization of compact and weakly compact intrinsic operators mapping into the growth spaces.

\bigskip

\textbf{Theorem 3.1.}\textit{ Let $T$ be an intrinsic operator on $\H(\D)$ and let $X$ be an initial space.} 

\textit{(i) Let $T: X \to H_{v}$ be bounded. Then $T: X \to H_v$ is compact (or weakly compact) if and only if the set 
$\{ v(z) T^* K_z^H; z \in \D \}$ is relatively compact (or relatively weakly compact) in $X^*$.}

\textit{(ii) Let $T: X \to \B_{v}$ be bounded. Then $T: X \to \B_{v}$  is compact (or weakly compact) if and only if the set 
$\{ v(z) T^* K_{z,1}^{\B}; z \in \D \}$ is relatively compact (or relatively weakly compact) in $X^*$.}

\begin{proof} 
\textit{(i)} One direction is trivial and follows by Schauder's theorem [8, p. 323] since $v(z)||K_z^H|| \le 1$, 
i.e. since the set 
$\{ v(z)K_z^H; z \in \D \}$ is a bounded set in $H_v^*$ and $T^*: H_v^* \to X^*$ is compact.

For the other direction, assume that $\{ v(z) T^* K_z^H; z \in \D \}$ is relatively compact in $X^*$. Suppose that $T: X \to H_v$ is not compact, i.e. by Lemma 3.1, suppose that $\exists \e_0>0$ and $\{f_n\}$ in $X$ with $||f_n||_X \le 1$, such that 
$$f_n \to 0 ~\mbox{in}~ \tau_{uc}, ~\mbox{as}~ n \to \infty, ~\mbox{but} ~  ||Tf_n||_v \ge  2\e_0, \forall n \in \N.$$
Thus $\exists \{z_n\} \subset \D$ such that $v(z_n)|Tf_n(z_n)| \ge  \e_0$, i.e. 
$v(z_n)|T^* K_{z_n}^H (f_n)| \ge  \e_0, \forall n \in \N.$

Since $\{ v(z) T^* K_z^H; z \in \D \}$ is relatively compact in $X^*$, there exists a subsequence $\{z_{n_k}\}$ and $m \in X^*$ such that 
as $k \to \infty$,
$$\sup\{|(v(z_{n_k})T^*K_{z_{n_k}}-m)(f)|; ||f||_X \le 1 \} \to 0.$$
Let $k_0 \in \N$ be such that $\forall k \ge k_0, \forall n \in \N$
$$|(v(z_{n_k})T^*K_{z_{n_k}} - m)(f_n)| = |(v(z_{n_k})Tf_n(z_{n_k}) - <f_n, m>| < \frac{\e_0}{4}.$$
Since $f_n \to 0$ in $\tau_{uc}$ and $T$ is intrinsic, $|Tf_n(z_{n_{k_0}})| \to 0$, as $n \to \infty$. Thus, for sufficiently large $n$
$$|<f_n,m>| \le |(v(z_{n_{k_0}})Tf_n(z_{n_{k_0}}) - <f_n, m>| + |(v(z_{n_{k_0}})Tf_n(z_{n_{k_0}})|< \frac{\e_0}{4}+\frac{\e_0}{4}=\frac{\e_0}{2},$$
and so $\limsup_{n \to \infty}|<f_n,m>| \le \frac{\e_0}{2}$. But then for large enough $k \ge k_0$,
\begin{eqnarray*}
\frac{5\e_0}{8} &>& |<f_{n_k}, m>| \ge |v(z_{n_k})Tf_{n_k}(z_{n_k})| - |v(z_{n_k})Tf_{n_k}(z_{n_k}) - <f_{n_k}, m>| \\
 &\ge& \e_0 -\frac{\e_0}{4}=\frac{6\e_0}{8} 
\end{eqnarray*}
and we get a contradiction.
Hence, by Lemma 3.1., $T:X \to H_v$ is compact.

The proof of \textit{(ii)} follows similarly. The only difference is that we have to estimate $|(Tf_n)'(z)|$ instead of $|(Tf_n)(z)|$. But since 
$(Tf_n)' \to 0$ in $\tau_{uc}$, whenever $Tf_n \to 0$ in $\tau_{uc}$, the rest of the proof follows in the same fashion.

\end{proof}

\smallskip

A similar type of characterization of compactness of general (intrinsic) operators on the so called admissible spaces, but with the disc algebra $A(\D)$ as the target space, was obtained in [5, Lemma 4.1]. Even though the methods and the important duality relations are different, in principle, the disc algebra role in [5] corresponds in many instances to the role played here by the little growth spaces. Contrary to that general principle, our next result gives a much more specific description of the compactness of intrinsic operators on general initial spaces, mapping into the little growth spaces, than the one given in [5] for the case of the disc algebra. 

\bigskip

\textbf{Theorem 3.2.} \textit{Let $T$ be an intrinsic operator on $\H(\D)$ and let $X$ be an initial space.}

\noindent \textit{(a) If  $T: X \to H_{v, 0}$ is bounded, then} 

\textit{(i) $T: X \to H_{v,0}$ is compact if and only if $\lim_{|z| \to 1} v(z) ||T^*K_z^H|| = 0$.}

\textit{(ii) $T: X \to H_{v, 0}$ is weakly compact if and only if $v(z) T^* K_z^H \to 0$ weakly,} 

\textit{\quad as $|z| \to 1$.}

\noindent \textit{(b) If $T: X \to \B_{v,0}$ is bounded, then} 

\textit{(i) $T: X \to \B_{v,0}$ is compact if and only if $\lim_{|z| \to 1} v(z) ||T^*K_{z,1}^{\B}|| = 0$.}

\textit{(ii) $T: X \to \B_{v,0}$ is weakly compact if and only if $v(z) T^* K_{z,1}^{\B} \to 0$ weakly,} 

\textit{\quad as $|z| \to 1$.}

\begin{proof}

We will show the proof only for case (a),  since case (b) follows similarly.

Part (i):
Since $T$ is intrinsic, $T: X \to H_{v,0}$ is bounded and $H_{v,0}$ is a closed subspace of $H_v$, $T: X \to H_{v,0}$ is compact if and only if $T: X \to H_{v}$ is compact, which by Theorem 3.1. is further more equivalent to $\{ v(z) T^* K_z^H; z \in \D \}$ is relatively compact in $X^*$. 

Hence, let us show first that $\lim_{|z| \to 1} v(z) ||T^*K_z^H|| = 0$ implies that $\{ v(z) T^* K_z^H; z \in \D \}$ is relatively compact in 
$X^*$, i.e. that for any sequence $\{z_n\} \subset \D$, there is a subsequence $\{z_{n_k}\}$ and $m \in X^*$ such that 
$v(z_{n_k}) T^* K_{z_{n_k}}^H \to m$ in $X^*$.

Case 1. If $\{z_n\}$ has a subsequence $\{z_{n_k}\}$ such that $z_{n_k} \to z_0 \in \D$, then 
$$v(z_{n_k}) T^* K_{z_{n_k}}^H \to v(z_0) T^*K_{z_0}^H = m \in X^*,$$
since $v, T^*$ and the map $z \to K_z^H$ are all continuous.

Case 2. If on the other hand $\{z_n\}$ is such that $|z_n| \to 1$ as $n \to \infty$, then
$$0 \le \lim_{n \to \infty} v(z_n) ||T^*K_{z_n}^H|| \le \lim_{|z| \to 1} v(z) ||T^*K_z^H|| = 0,$$
and so $v(z_{n}) T^* K_{z_{n}}^H \to 0$ in $X^*$.

For the other direction, assume that $\{ v(z) T^* K_z^H; z \in \D \}$ is relatively compact in $X^*$. Suppose that there exists a sequence
$\{z_n\} \subset \D$ with $|z_n| \to 1$ as $n \to \infty$, and $\exists c>0$ such that $v(z_n) ||T^*K_{z_n}^H|| \ge c>0$. Choose a sequence
$\{f_n\} \subset B_X$ such that
$$v(z_n) |T^*K_{z_n}^H(f_n)| = v(z_n) |Tf_n (z_n)| > \frac{c}{2}.$$
Since $B_X$ is compact with respect to $\tau_{uc}$, without loss of generality, $f_n \to f$ in $\tau_{uc}$, for some $f \in X$.
But $T:X \to H_v$ is compact, and so by Lemma 3.3, $||T(f_n - f)||_v \to 0$. Since also $T:X \to H_{v,0}$ is bounded, $Tf_n \in H_{v,0}$,
and so it must be that $Tf \in H_{v,0}$. Furthermore, for sufficiently large n,
$$v(z_n) |Tf(z_n)| \ge v(z_n) |Tf_n (z_n)| - v(z_n) |T(f_n-f) (z_n)| \ge \frac{c}{2}-\frac{c}{4}=\frac{c}{4}>0.$$
This contradicts $Tf \in H_{v,0}$ since $\{z_n\}$ is such that  $|z_n| \to 1$ as $n \to \infty$.

Thus we must have that $\lim_{|z| \to 1} v(z) ||T^*K_z^H|| = 0$.

\medskip

Part (ii): As in the proof of part (i), since $T: X \to H_{v, 0}$ is bounded, $T: X \to H_{v, 0}$ is weakly compact if and only
if $T: X \to H_{v}$ is weakly compact.
 
So, let $v(z) T^* K_z^H \to 0$ weakly, as $|z| \to 1$. By Theorem 3.1, we need to show that then the set $\{ v(z) T^* K_z^H; z \in \D \}$ is relatively weakly compact in $X^*$. By Eberlein-$\check{\text S}$mulian theorem [13, p. 49], we have to show that every sequence in this set has a weakly convergent subsequence, i.e.
 that for any sequence $\{z_n\} \subset \D$, there is a subsequence $\{z_{n_k}\}$ and $m \in X^*$ such that 
 $v(z_{n_k}) T^* K_{z_{n_k}}^{H} \to m$ weakly in $X^*$. As in part (i), there are two possible cases to consider:

Case 1. If $\{z_n\}$ has a subsequence $\{z_{n_k}\}$ such that $z_{n_k} \to z_0 \in \D$, then 
$$v(z_{n_k}) T^* K_{z_{n_k}}^H \to v(z_0) T^*K_{z_0}^H = m \in X^*,$$
since $v, T^*$ and the map $z \to K_z^H$ are all continuous. But then 
$v(z_{n_k}) T^* K_{z_{n_k},1}^{\B} \to m$ also weakly in $X^*$.

Case 2. If on the other hand $\{z_n\}$ is such that $|z_n| \to 1$ as $n \to \infty$, then trivially, by our assumption, 
$v(z_{n_k}) T^* K_{z_{n_k}}^{H} \to 0$ weakly in $X^*$.

For the other direction, let $T: X \to H_{v, 0}$ be weakly compact.
By Gantmacher's Theorem [8, p.343], $T^*:H_{v,0}^* \to X^*$ is also weakly compact. 
So, we need to show that $T^*:H_{v,0}^* \to X^*$ weakly compact implies that $v(z)T^*K_z^H$ converges weakly to $0$, as $|z| \to 1$.
Since $||K_z^H|| \le 1/v(z)$, i.e. $v(z)K_z^H$ is in the closed unit ball of $H_{v,0}^*$, the set $\{v(z)T^*K_z^H; z \in \D \}$ is relatively weakly compact in $X^*$. Hence, by Eberlein-$\check{\text S}$mulian theorem, every sequence in this set has a weakly convergent subsequence.

Suppose $v(z)T^*K_z^H$ does not converge weakly to $0$ in $X^*$, as $|z| \to 1$, i.e. there exists $F \in X^{**}$ such that
the (bounded) set $\{|F(v(z)T^*K_{z}^H )|; z \in \D \}$ does not converge to $0$ in $\R^{+}$, as $|z| \to 1$. 
Let $\{z_n\}$ in $\D$, $|z_n| \to 1$, be such that  
$$|F(v(z_n)T^*K_{z_n}^H )| \geq c > 0.$$ 
From the relative weak compactness, there exists a subsequence $\{z_{n_k}\}$, such that  $\{v({z_{n_k}})T^*K_{z_{n_k}}^H \}$ converges weakly to some $m \in X^*$. 
But then the last subsequence also converges weak$^*$ to $m$, i.e. $\forall f \in X$ 
$$\lim_{k \to \infty} v({z_{n_k}})T^*K_{z_{n_k}}^H(f)=\lim_{l \to \infty} v({z_{n_k}})Tf(z_{n_k})=m(f).$$
Since $T$ is bounded, $Tf \in H_{v,0}$, i.e. 
$$\lim_{k \to \infty} v({z_{n_k}})|Tf(z_{n_k})|=0=|m(f)|,$$ 
for all $f \in X$. Thus, $m=0$ and so $\lim_{k \to \infty}F(v({z_{n_k}})T^*K_{z_{n_k}}^H)=0$, which contradicts 
$|F(v(z_n)T^*K_{z_n}^H )| \geq c > 0$.

\end{proof}

\medskip

The criteria in the previous result have a familiar form when applied to some specific operators. This is true in particular for the compactness case, and we state the special case of weighted composition operators and integral operators as an illustration of the application of Theorem 3.2.
One of the first results of this kind was the result of Madigan and Matheson in [7] for the composition operators on the little Bloch space. A more general version for weighted composition operators on initial spaces with bounded point evaluations, and mapping into the spaces $H_{v,0}$, namely part (i) of (a) in the following corollary, appears in [4]. 

\bigskip

\textbf{Corollary 3.1.}\textit{ Let $X$ be an initial space with bounded point evaluations $\delta_z$, and let $\b>0$.}

\noindent \textit{(a) If the weighted composition operator $W_{u,\p}: X \to H_{v,0}$ is bounded, then} 

\textit{(i) $W_{u,\p}: X \to H_{v,0}$ is compact if and only if $\lim_{|z| \to 1} v(z) |u(z)| ||\delta_{\p(z)}|| = 0$.}

\textit{(ii) $W_{u,\p}: X \to H_{v,0}$ is weakly compact if and only if $v(z) u(z) \delta_{\p(z)} \to 0$ weakly,} 

\textit{\quad as $|z| \to 1$.}

\noindent \textit{(b) If the integral operator $T_g: X \to \B_{v,0}$ is bounded, then} 

\textit{(i) $T_g: X \to \B_{v,0}$ is compact if and only if $\lim_{|z| \to 1} v(z) |g'(z)| ||\delta_z|| = 0$.}

\textit{(ii) $T_g: X \to \B_{v,0}$ is weakly compact if and only if $v(z) g'(z) \delta_z \to 0$ weakly,} 

\textit{\quad as $|z| \to 1$.}

\noindent \textit{(c) If the integral operator $T_g: X \to H_{\b, 0}$ is bounded, then} 

\textit{(i) $T_g: X \to H_{\b,0}$ is compact if and only if $\lim_{|z| \to 1} (1-|z|^2)^{\b+1} |g'(z)| ||\delta_z|| = 0$.}

\textit{(ii) $T_g: X \to \B_{v,0}$ is weakly compact if and only if $(1-|z|^2)^{\b+1} g'(z) \delta_z \to 0$} 

\textit{\quad weakly, as $|z| \to 1$.}

\bigskip

To address further the compactness and weak compactness of intrinsic operators $T: X \to H_v (B_v)$, we look at a special class of initial spaces. The extra conditions on the space $X$ imposed in the next few theorems are part of the restrictions on the "admissible" domains considered in [4] and [5] for the weighted composition operators and the integral operators. 
We start first by considering a special kind of intrinsic weakly compact operators. Note that a similar result ([5, Proposition 3.4]) holds for the integral operators, with $H^{\infty}$ and $A(\D)$ playing the roles of $H_v (B_v)$ and $H_{v,0}$ below.

\bigskip

\textbf{Theorem 3.3.}\textit{ Let $X$ be an initial space such that for $f \in X$ and $0 \leq r<1$ the functions 
$f_r(z)=f(rz)$ are in $X$, and furthermore
$$\sup_{0 \leq r<1} ||f_r||_X \lesssim ||f||_X.$$}
\textit{ Let $T$ be an intrinsic operator on $\H(\D)$ satisfying $Tf_r \in H_{v,0}$, for all $f \in X$ and all $0 \le r <1$.
Then $T: X \to H_{v}$ weakly compact implies that $T: X \to H_{v,0}$ is bounded.} 

\begin{proof} By our assumptions, for $f \in X$ and $0 \le r <1$, we have that $f_r \in X$. Also, in  general, $f_r$ converges to $f$ in $\tau_{uc}$, as $r \to 1$.
Thus, since $T: X \to H_v$ is weakly compact, we have from Lemma 3.1 that 
$Tf_r$ converges weakly in $H_v$ to $Tf$. Since $H_{v,0}$ is a closed subspace of $H_v$, by Mazur's theorem ([13, p. 28]) the weak and the norm closure of $H_{v,0}$ coincide. But $Tf_r \in H_{v, 0}$, $Tf_r$ converges weakly in $H_v$ to $Tf$, and so
$Tf \in H_{v,0}$.  Thus,  $T(X) \subset H_{v,0}$ and by the closed graph theorem, $T:X \to H_{v,0}$ is bounded.

\end{proof}

Note that since $v$ is a typical weight, for any $f$ in $X$ each of the functions $f_r$ is actually contained on $H_{v,0}$. Requiring that 
$H_{v,0} \subset X$, and that $H_{v,0}$ is invariant for $T$ is also a natural possible condition, but that is in general a stronger condition than the one given in the previous theorem. 

\bigskip

A similar results hold for the case of Bloch type growth spaces. Since the proof is very similar, we will only state the corresponding results and leave the proof to the reader.

\bigskip

\textbf{Theorem 3.3'.} \textit{ Let $X$ be an initial space such that for $f \in X$ and $0 \leq r<1$ the functions 
$f_r(z)=f(rz)$ are in $X$, and furthermore
$$\sup_{0 \leq r<1} ||f_r||_X \lesssim ||f||_X.$$}
\textit{ Let $T$ be an intrinsic operator on $\H(\D)$ satisfying $Tf_r \in \B_{v,0}$, for all $f \in X$ and all $0 \le r <1$.
Then $T: X \to \B_{v}$ weakly compact implies that $T: X \to \B_{v,0}$ is bounded.}

\bigskip

In the previous two theorems, since the space $X$ contains the constants, the condition $Tf_r \in H_{v,0}(\B_{v,0})$ also implies that 
$T1 \in H_{v,0}(\B_{v,0})$. For the case of the weighted composition operator $W_{u,\p}$, this means that $u \in H_{v,0}(\B_{v,0})$, and if
$T$ is the integral operator $T_g$, we get that $g \in H_{v,0}(\B_{v,0})$. Note that both are necessary conditions for the boundedness of $T: X \in H_{v,0}(\B_{v,0})$. On the other hand, it is not hard to see that if either $u \in H_{v,0}$, or $g \in \B_{v,0}$, then $W_{u,\p}$, or $T_g$ correspondingly, satisfy the required conditions in Theorem 3.3 and Theorem 3.3', namely that then $W_{u,\p} f_r \in H_{v,0}$, or
$T_g f_r \in \B_{v,0}$, for all $f \in X$ and all $0 \le r <1$. 

\bigskip

Using the fact that every bounded operator on a reflexive space is also weakly compact and that $H_{v,0}^{**} =H_v$ and  
$\B_{v,0}^{**} =\B_v$, we get the following corollary. 

\bigskip

\textbf{Corollary 3.2.}\textit{ Let $X$ and $T$ satisfy the conditions of Theorem 3.3 (or Theorem 3.3'). If $X$ is also a reflexive space, then $T: X \to H_{v}(\B_v)$ is bounded if and only if $T: X \to H_{v,0}(\B_{v,0})$ is bounded.}
\begin{proof} 
If $T: X \to H_{v,0}$ is bounded, since $X=X^{**}, T=T^{**}$ and $H_{v,0}^{**} =H_v$, we get that also $T: X \to H_{v}$ is bounded.
On the other hand, if $T: X \to H_{v}$ is bounded, it is also weakly compact, since $X$ is reflexive. But then by Theorem 3.3, 
$T: X \to H_{v,0}$ is bounded.

The proof of the case $\B_v$ and $\B_{v,0}$ follows similarly.

\end{proof} 

\bigskip

Note that a similar result holds in [5] for integral operators mapping a reflexive, admissible space $X$ 
into $H^{\infty}$. It is interesting that in that case, as before, the space $H_{v, 0}$ is replaced by the disk algebra $A(\D)$, even though the double dual of  $A(\D)$ only includes $H^{\infty}$, and is not equal to  $H^{\infty}$. (See Theorem 1.1, part (ii) and Corollary 3.8, part (ii) in [5].)

\bigskip

Next, we would like to get a more specific description of compact intrinsic operators mapping into the spaces $H_v$ and $\B_v$, when $X$ is not necessarily reflexive. We have the following sufficient condition.

\bigskip

\textbf{Proposition 3.1.} \textit{Let $T$ be an intrinsic operator on $\H(\D)$ and let $X$ be an initial space.}
 
\textit{(i) Let $T: X \to H_{v}$ be bounded. For $N \in \N$, let $D_N=\{z \in \D; ||T^*K_z^H||>N\}$. 
Then 
\begin{equation*}
\lim_{N \to \infty} \sup_{z \in D_N}v(z) ||T^*K_z^H|| = 0 
\tag{3.1.1}
\end{equation*}
implies that $T: X \to H_{v}$ is compact.}

\textit{(ii) Let $T: X \to \B_{v}$ be bounded. For $N \in \N$, let $D_N=\{z \in \D; ||T^*K_{z,1}^{\B}||>N\}$. 
Then 
\begin{equation*}
\lim_{N \to \infty} \sup_{z \in D_N}v(z) ||T^*K_{z,1}^{\B}|| = 0 
\tag{3.1.2}
\end{equation*}
implies that $T: X \to \B_{v}$ is compact.}

\begin{proof} Part (i): 
We will show that (3.1.1) implies that the set $\{ v(z) T^* K_z^H; z \in \D \}$ is relatively compact in $X^*$. Then the compactness of  $T: X \to H_{v}$ follows from Theorem 3.1. 

As in the proof of the previous theorem, start with a general sequence 
$\{z_n\} \subset \D$. 
If $\{z_n\}$ has a subsequence $\{z_{n_k}\}$ such that $z_{n_k} \to z_0 \in \D$, then using the continuity of $v, T^*$ and the map 
$z \to K_z^H$ we get that
$$v(z_{n_k}) T^* K_{z_{n_k}}^H \to v(z_0) T^*K_{z_0}^H = g \in X^*.$$

Otherwise, $|z_n| \to 1$, as $n \to \infty$. In case $\exists c>0$ and a subsequence $\{z_{n_k}\}$ with 
$||T^*K_{z_{n_k}}|| \le c$, using that $v$ is a typical weight, we have that $v(z_{n_k}) T^* K_{z_{n_k}}^H \to 0$ in $X^*$.
On the other hand if $||T^*K_{z_{n}}|| \to \infty$, without loss of generality, pick a subsequence $\{z_{n_N}\}$ such that 
$z_{n_N} \in D_N$. But then
$$0 \le \lim_{N \to \infty} v(z_{n_N}) ||T^*K_{z_{n_N}}^H|| \le \lim_{N \to \infty} \sup_{z \in D_N}v(z) ||T^*K_z^H|| = 0,$$
and so $v(z_{n_N}) T^* K_{z_{n_N}}^H \to 0$ in $X^*$. Thus, $\{ v(z) T^* K_z^H; z \in \D \}$ is relatively compact in $X^*$, and
$T: X \to H_{v}$ is compact.

The proof of part (ii) follows similarly, this time using the derivative point evaluations  $K_{z,1}^{\B}$ instead of the point evaluations $K_z^H$.

\end{proof}

\bigskip

It would be interesting to know if in general, the condition (3.1.1) (or (3.1.2)) is also a necessary condition for compactness of intrinsic operators. Proving the more specific cases of this problem usually requires a proper choice of test functions that depend on the operator and/or on the space $X$. Since the class of initial spaces and the class of intrinsic operators are quite general classes, determining this might be out of reach. 
On the other hand, there is more hope if one is interested in a specific class of operators. For example, the results for the class of weighted composition operators in [4] show that if $X$ is an initial space with bounded point evaluations, satisfying few more specific properties, then the condition (3.1.1) (or (3.1.2)) is also a necessary compactness condition for bounded weighted composition operators (see Theorem 3.2 and Theorem 3.3 in [4]).

\bigskip

The next result shows that (3.1.1) and (3.1.2) are necessary and sufficient conditions for compactness also of nicely behaved intrinsic operators, acting on a smaller class of "admissible" initial spaces $X$, namely for the spaces $X$ also considered in Theorem 3.3 and Theorem 3.3' above.

\bigskip

\textbf{Theorem 3.4.}\textit{ Let $X$ be an initial space such that for $f \in X$ and $0 \leq r<1$ the functions 
$f_r(z)=f(rz)$ are in $X$, and furthermore
$\sup_{0 \leq r<1} ||f_r||_X \lesssim ||f||_X$. Let $T$ be an intrinsic operator on $\H(\D)$, $T: X \to H_v$ bounded, and $Tf_r \in H_{v,0}$, for all $f \in X, 0 \le r <1$. Then the following are equivalent:}

\textit{(i) $T: X \to H_{v}$ is compact.} 
 
\textit{(ii) $T: X \to H_{v,0}$ is compact.}
 
 \textit{(iii) $\lim_{|z| \to 1} v(z) ||T^*K_z^H|| = 0$.}
 
 \begin{proof} We will show that $(i)\Rightarrow(ii)\Rightarrow(iii)\Rightarrow(i)$.
 
 If $T: X \to H_{v}$ is compact, then $T: X \to H_{v}$ is weakly compact and by Theorem 3.3, $T: X \to H_{v,0}$ is bounded. Since 
 $H_{v,0}$ is a closed subset of $H_v$, $T: X \to H_{v, 0}$ is compact. Then, by Theorem 3.2, $\lim_{|z| \to 1} v(z) ||T^*K_z^H|| = 0$. It is easy to see that this implies $\lim_{N \to \infty} \sup_{z \in D_N}v(z) ||T^*K_z^H|| = 0$, which by Proposition 3.1 implies that $T: X \to H_{v}$ is compact.
 
 \end{proof}
 
 \bigskip
 
 As before, a similar result also holds for the Bloch type growth spaces. 
 
 \bigskip
 
 \textbf{Theorem 3.4'.}\textit{ Let $X$ be an initial space such that for $f \in X$ and $0 \leq r<1$ the functions 
$f_r(z)=f(rz)$ are in $X$, and furthermore
$\sup_{0 \leq r<1} ||f_r||_X \lesssim ||f||_X$. Let $T$ be an intrinsic operator on $\H(\D)$, $T: X \to \B_v$ bounded, and $Tf_r \in \B_{v,0}$, for all $f \in X, 0 \le r <1$. Then the following are equivalent:}

\textit{(i) $T: X \to \B_{v}$ is compact.} 
 
\textit{(ii) $T: X \to \B_{v,0}$ is compact.}
 
 \textit{(iii) $\lim_{|z| \to 1} v(z) ||T^*K_{z,1}^{\B}|| = 0$.}
 
  \bigskip
  
 We illustrate an application of Theorem 3.4 for the integral operator $T_g$. Note that when $X$ contains the constants, then $T_g f_r \in H_{\b,0}$ for all $f \in X, 0 \le r <1$ if and only if $g \in H_{\b,0}$. To show also the connections with some of the boundedness criteria in Corollary 2.2, and to further flag the similarities with the compactness results in [5, Theorem 4.2 and Theorem 4.3] for the integral operators mapping into $H^{\infty}$ and $A(\D)$ spaces, we state just the following two particular results.
 
 \bigskip

\textbf{Corollary 3.3.}\textit{ Let $1<p< \infty$, $\a>-1$, and $\b>0$.}

\noindent \textit{(a) If the integral operator $T_g: H^p \to H_{\b}$ is bounded and $g \in H_{\b,0}$, then the following are equivalent:} 

\textit{(i) $T_g: H^p \to H_{\b}$ is compact.}

\textit{(ii) $T_g: H^p \to H_{\b,0}$ is compact.} 

\textit{(iii) $\lim_{|z| \to 0}(1-|z|^2)^{\beta+1 - \frac{1}{p}} |g'(z)|=0$.}

\noindent \textit{(b) If the integral operator $T_g: A_{\a}^p \to H_{\b}$ is bounded and $g \in H_{\b,0}$, then the following are equivalent:} 

\textit{(i) $T_g: A_{\a}^p \to H_{\b}$ is compact.}

\textit{(ii) $T_g: A_{\a}^p \to H_{\b,0}$ is compact.} 

\textit{(iii) $\lim_{|z| \to 0}(1-|z|^2)^{\beta+1 - \frac{\a+2}{p}} |g'(z)|=0$.}

\bigskip

Notice, for example, that we can conclude from part (b), (iii), that when $p =\frac{\a+2}{\b+1}$, $T_g: A_{\a}^p \to H_{\b}$ 
(or $T_g: A_{\a}^p \to H_{\b,0}$) is compact only when $g$ is a constant, while, for example, $T_g$ is bounded for any polynomial $g$. 
(See [5] for further similar conclusions, and other such particular situations.)

If in Theorem 3.4 we take $T$ to be a weighted composition operator $W_{u,\p}$, and $X=H_w$, we get the results that were previously derived in [3] (see Theorem 3.3) for the compactness of composition operators $C_{\p}: H_w \to H_v$, or in [9] (see Theorem 2.2) for the weighted composition operators $W_{u,\p}: H_w \to H_v$.
  
\bigskip
 
Another interesting problem is to determine when are the weak compactness and compactness equivalent for a given operator. This is definitely the case if the domain of the operator is a "$c_0$ type" Banach space (see [8, p.347]). Such are, for example the spaces $H_{v,0}$ and $\B_{v,0}$. Note that they are not initial spaces. Still, it might be possible that under some special conditions on the initial space $X$ and on the intrinsic operator $T$, we have the equivalence of the compactness and weak compactness of $T$, when mapping such $X$ into the growth, or little growth spaces. We leave this as an open question for further consideration.

\bigskip

\bigskip


\medskip
 
\begin{thebibliography}{91}

\bibitem{AlPe} A. Aleman, A. Persson, \textit{Resolvent estimates and decomposable extensions of generalized Ces\'aro operators}, J. Funct. Anal. 258 (2010), 67-98.

\bibitem{Bie} K. D. Bierstedt, J. Bonet, J. Taskinen, \textit{Associated weights and spaces of holomorphic functions}, Studia Math.
127 (1998), 137-168.

\bibitem{Bon} J. Bonet, P. Domanski, M. Lindstrom, J. Taskinen, \textit{Composition operators between weighted Banach spaces of analytic functions}, J. Austral. Math. Soc. 64 (1998), 101-118.

\bibitem{CoTj} F. Colonna, M. Tjani, \textit{Operator norms and essential norms of weighted composition operators between Banach spaces of analytic functions}, J. Math. Anal. Appl. 434 (2016), 93-124.

\bibitem{CoPel} M. D. Contreras, J. A. Pel\'aez, C. Pommerenke, J. Rattya, \textit{Integral operators mapping into the space of bounded analytic functions}, J. Funct. Anal. 271 (2016), 2899-2943.

\bibitem{Hu} Z. Hu, \textit{Extended Cesaro operators on mixed norm spaces}, Proc. Amer. Math. Soc., Vol. 131, No. 7 (2003), 
2171-2179.

\bibitem{MadMa} K. Madigan, A. Matheson \textit{Compact composition operators on the Bloch space}, Trans. Amer. Math. Soc., Vol.347, No. 7 (1995), 2679- 2687.

\bibitem{Megg} R. E. Megginson, An Introduction to Banach Space Theory, Graduate Texts in Mathematics, Vol. 183, Springer, 1998.

\bibitem{Mo} A. Montes-Rodriguez, \textit{Weighted composition operators on weighted Banach spaces of analytic functions}, J. London Math. Soc. (2) 61 (2000), 872-884.

\bibitem{OhStr} S. Ohno, K. Stroethoff, \textit{Weighted composition operators from reproducing Hilbert spaces to Bloch spaces}, Houston J. Math., Vol. 37, No. 2 (2011), 537-558.

\bibitem{Ste} S. Stevi\'c, \textit{On an integral operator on the unit ball in $\C^n$}, J. Inequal. Appl. 2005 (2005), 81-88.

\bibitem{Tj} M. Tjani, \textit{Compact composition operators on Besov spaces}, Trans. Amer. Math. Soc., Vol. 355, No. 11 (2003), 
4683-4698.

\bibitem{Wojt} P. Wojtaszczyk, Banach Spaces for Analysts, Cambridge Studies in Advanced Mathematics, vol. 25,  Cambridge University Press, Cambridge, 1991.

\bibitem{Zhu} K. Zhu, \textit{Bloch type spaces of analytic functions}, Rocky Mountain J. Math. 23 (3) (1993), 1143-1177.

\end{thebibliography}
\end{document}